\documentclass[12pt]{amsart}
\usepackage{geometry}          % See geometry.pdf to learn the layout options.
\usepackage{graphicx}
\usepackage{amssymb}
\usepackage{epstopdf}
\def\F{{\Bbb{F}}}

\def\N{{\Bbb N}}

\def\TT{{\bf T}}\def\TF{{\TT^F}}\def\SS{{\bf S}}

\def\GG{{\bf G}}\def\GF{{{\bf G}^F}}

\def\LL{{\bf L}}

\def\aa{{\alpha}}

\def\dd{{\delta}}

\def\si{{\sigma}}

\def\th{{\theta}}
\def\la{{\lambda}}
\def\La{{\Lambda}}

\def\Aut{{\rm Aut}}

\def\Ind{{\rm Ind}}
\def\Ce#1{{{\rm C}_#1 }}

\def\No#1{{{\rm N}_{#1}}}

\def\inn{\subseteq}

\def\vid{\emptyset}

\def\ser#1#2{\mathcal{E} (#1 ,#2 )}

\def\Sp{{\rm Sp}}
\def\GL{{\rm GL}}
\def\Irr{{\rm Irr}}
\def\lp{{\ell '}}

\def\inn{{\subseteq}}

\def\Lu#1#2{{\rm R}_{#1}^{#2}}
\def\qed{\hfill\vrule height 1.6ex width .7ex depth -.1ex }

\def\mr#1{\smash{\mathop{\relbar\joinrel\longrightarrow}\limits^{#1}}}

\title{Odd character degrees for $\Sp (2n,2)$.}
\author{Marc Cabanes}

\address{Institut de Math\'ematiques de Jussieu, Universit\'e Paris 7, 175 rue du Chevaleret,  F-75013 Paris, FRANCE}

\email{cabanes@math.jussieu.fr}

\begin{document}

\begin{abstract}
We check McKay conjecture on character degrees for the case of symplectic groups over the field with two elements $\Sp_{2n}(2)$ and the prime 2. Then we check the inductive McKay condition (see [IMN] 10) for $\Sp_4(2^m)$ and all primes.
\end{abstract}

\maketitle

%\centerline{}
%\bigskip
%{}
 If $G$ is a finite group and $\ell$ is a prime number, denote by $\Irr_{\ell '}(G)$ the set of irreducible characters of $G$ with
degree prime to $\ell$. The McKay conjecture on
character degrees asserts that $$|\Irr_{\ell '}(G)|=|\Irr_{\ell '}(\No G(P))|$$ for $P$ a Sylow
$\ell$-subgroup of $G$. McKay's conjecture has gained new interest
since appearance of Isaacs-Malle-Navarro's theorem reducing it to a related
conjecture on quasi-simple groups
(see [IMN]).  It has been checked for all quasi-simple
groups not of Lie type.

Among groups of Lie type and for the prime $\ell$ being the defining characteristic, the
group $\Sp_{2n}(2)$ is one of the rare cases, and the only infinite
series, of finite groups of Lie type whose Sylow subgroup at
defining characteristic (2 in the present case) has an abelian quotient bigger than
expected from the rank of the underlying reductive group (see
Proposition~3 below). We prove that it nevertheless satisfies McKay conjecture, which in this case also finishes the checking of the conditions devised by [IMN]. For a uniform treatment of the ``general" case, see [S3].

We intended to give a sequel to this note by checking the conditions of [IMN] for all simple groups Sp$_{2n}(2^m)$ but the equivariance problems proved more delicate than we first thought. In a joint work with B. Sp\"ath, we developed some general methods which cover this case (see [CS]). 
Here, we present however the case of Sp$_4(2^m)$ which requires some ad hoc analysis (see \S 2).

\section{Notations}

When $\ell$ is a prime and $n\geq 1$ an integer, one denotes by $n_\ell$ the greatest power of $\ell$ dividing $n$ and $n_\lp:=n/n_\ell$. If $H$ is a finite group and $X\inn \Irr (H)$, one denotes $X_\lp :=X\cap \Irr_\lp (H)$. 

If $H$ acts on a set $Y$, one denotes by $Y^H$ the subset of fixed points. If $H$ is cyclic generated by $h$, one writes $Y^H=Y^h$. If $Y$ is a group on which $H$ acts by group automorphisms, one denotes by $[Y,H]$ the subgroup of $Y$ generated by commutators $y^{-1}h(y)$ for $h\in H$, $y\in Y$.

For finite reductive groups $\GF$ ($\GG$ a reductive group defined over a finite field $\F_q$ with associated Frobenius endomorphism $F\colon \GG\to\GG$), the Deligne-Lusztig generalized characters $\Lu\TT\GG\th$, and associated partition of $\Irr (\GF)$ into series $\ser\GF s$ (in our cases rational and geometric series will always coincide since $\GG$ will have connected center) we follow the notations of [DM].

\section{Odd character degrees for $\Sp_{2n}(2)$}

Let us denote by $\F$ the algebraic
closure of $\F_2$ the field with 2 elements. Let $n\geq 2$ be an integer,
let $\GG =\Sp_{2n}(\F )$ with Frobenius endomorphism ${F_0}\colon
\GG\to \GG$ squaring matrix entries. Let
$G=\GG^{F_0}=\Sp_{2n}(\F_2)$.

\medskip

\subsection{The global case}

We use the usual notion and notation $\mathcal{E}(G,1)$ for unipotent
characters
(see [L], [C2]~\S~12, [DM]~13.19). Unipotent characters of finite
classical groups over a finite field are described by the work of
Lusztig [L] (see also [C2]\S 13.8). They are parametrized by a set
(independent of the finite field) of so-called symbols
$\La$, which in the case of type
$B_n$ are of the kind $\La = \left(\begin{matrix}  S\\ T \end{matrix} \right)$  where
$S,T$ are subsets of $\N$ with $0\not\in S\cap T$, $|S|-|T|$ is odd
and positive, and $n=\sum_{\la\in S}\la+\sum_{\mu\in T}\mu
-({|S|+|T|-1\over 2})^2$. Denote by $\chi_\La\in \mathcal{E}(G,1)$ the
unipotent character associated with $\La$. For the following, see also [M1]~6.8 but we give a full proof for the convenience of the reader.

\medskip\noindent {\bf Proposition~1.} {\sl  Keep
$n\geq 2$. Then $\chi_\Lambda (1)$ is odd if and only if
$\Lambda$ is among the {\bf five} following symbols $\left(\begin{matrix}
n
\\   \\ \end{matrix} \right)  $, $\left(\begin{matrix}  0,1,n \\
\\ \end{matrix} \right)  $,
$\left(\begin{matrix}  0,1 \\  n \\ \end{matrix} \right)  $,
$\left(\begin{matrix}  1,n \\  0 \\ \end{matrix} \right)  $, and $\left(\begin{matrix}
0,n \\  1
\\ \end{matrix} \right)  $ .  }

\medskip

\medskip\noindent {\it Proof.} For a unipotent character of a finite
classical group over a finite field of cardinality $q$ and
characteristic
$p$ (see a list of types in [L]~\S~8), $\chi_\La(1)_p = D_\La (q)_p$
where
$D_\La$ is the rational function in $q$ given in [L]~2.8 and
[C2]~p.467. In particular for $p=2$, denoting $\La  = \left(\begin{matrix}  S\\
T \end{matrix} \right)$ and
$M=|S|+|T|$, we have
$$\chi_\La (1)_2={\Pi_{\la '<\la\ \!
 {\rm in}\ \! S\ \! }(2^\la-2^{\la '})_2\ \!\Pi_{\mu '<\mu \ \!
{\rm in}\ \! T\ \! }(2^\mu-2^{\mu '})_2\ \!\Pi_{\la \in S,\mu\in
T}(2^\la+2^{\mu})_2\over 2^{{M-1\over 2}}2^{{1\over 2}\bigl(
(M-2)(M-3)+(M-4)(M-5)+\dots \bigr)}}.$$

The denominator is equal to $2^{{M-1\over
24}(2M^2-7M+15)}$ for any odd $M\geq 1$ (easy induction).

 Concerning the numerator, if
$c,d\geq 0$ are integers, one has
$(2^c-2^d)_2=(2^c+2^d)_2=2^{{\rm min}(c,d)}$ whenever
$c\not=d$, while $2^c+2^d=2^{{\rm min}(c,d)+1}$ whenever $c=d$.
So $$\chi_\La (1)_2=2^{\phi_{S,T}-{{M-1\over 24}(2M^2-7M+15)}}$$
where
$$\phi_{S,T}=|S\cap T|+\sum_{\la '<\la
\ \rm in\ S}{\rm min}(\la ,{\la '})+
\sum_{\mu '<\mu
\ \rm in\ T}{\rm min}(\mu ,{\mu '})+\sum_{\la\in S,\mu
\in T}{\rm min}(\la ,{\mu }).$$

It is now clear that $\chi_\La (1)$ is odd for $ \Lambda\in\Big\{
\left(\begin{matrix}
n
\\   \\ \end{matrix} \right)  ,  \left(\begin{matrix}  0,1,n \\
\\ \end{matrix} \right)   ,
 \left(\begin{matrix}  0,1 \\  n \\ \end{matrix} \right)   ,
 \left(\begin{matrix}  1,n \\  0 \\ \end{matrix} \right)   , \left(\begin{matrix}  0,n \\
 1
\\ \end{matrix} \right)   \Big\}$. Conversely, if
$M=3$ and
$\chi_\La(1)$ is odd, then
$\La = \left(\begin{matrix} \la '<\la\\  \mu \end{matrix} \right)$ and $\phi_\La=\la
'+{\rm min} (\mu,\la ')+{\rm min}(\mu ,\la)$, or $\La =
\left(\begin{matrix} \la ''<\la'<\la\\   \\ \end{matrix} \right)$ and $\phi_\La =
\la '+2\la ''$. Having $\phi_\La =1$ clearly forces $\{\la
',\mu\}=\{0,1\}$ (since $\la '=\mu =0$ is forbidden) and $(\la '',\la
')=(0,1)$, respectively. Hence the Proposition for $M\leq 3$.

There remains to check that $\phi_{S,T}>{{M-1\over
24}(2M^2-7M+15)}$ as soon as $M>3$. Let us re-index the elements
of $S$ and $T$ in a single sequence $\nu_1\leq
\nu_2\leq\dots\leq \nu_{M}
$ such that $S$ and $T$ correspond with disjoint subsets of indexes
in $[1,M]$. Then
$$\phi_{S,T}=|S\cap T|+ \sum_{i<j}{\rm min}(\nu_i,\nu_j)=|S\cap
T|+\sum_{i=1}^M\nu_i(M-i).$$ Since the sequence $(\nu_i)$ is the
merging of two strictly increasing sequences, it may take a given
value of $S\cup T$ only once or twice, and only once for the value 0.
So $\nu_i\geq$ the i-th term of the sequence $0,1,1,2,2,\dots
{M-1\over 2}, {M-1\over 2}$. Then $\phi_{S,T}\geq |S\cap
T|+\sum_{i=1}^{ {M-1\over 2}}i(2M        -4i-1)$. An easy induction on ${M-1\over 2}$
shows that
$$\sum_{i=1}^{ {M-1\over 2}}i(2M -4i-1) -
{M-1\over 24}(2M^2-7M+15)={(M-1)(M-3)\over 4}$$ for any odd $M\geq 3$. Thus
our claim.
\qed

\medskip

\noindent {\bf Proposition~2. } {\sl Let
$n\geq 2$ be an integer. Then $\Sp_{2n}(\F_2)$ has
$2^{n+1}$ characters of odd degrees.}

\medskip\noindent {\it Proof. } Recall $\GG =\Sp_{2n}(\F )$ with
Frobenius endomorphism ${F_0}\colon
\GG\to \GG$ squaring matrix entries. Let
$G=\GG^{F_0}=\Sp_{2n}(\F_2)$ (part of case (a) in [L]~\S~8). Note that
$\GG$ has (trivial) connected center.

By [L] p.164, $\F$
being of characteristic 2, there is an isogeny between $\GG$ and its
dual
$\GG^*$ inducing a bijection between rational semi-simple
elements with isomorphism of centralizers of corresponding elements. This, along with
property (A) of [L]~7.8 shows that ${\rm Irr} (G)$ is in bijection
with the disjoint union of the $\ser{\Ce G(s)}1$'s for $s$ ranging
over the semi-simple conjugacy classes of $G$ (see [L]~8.7.6).
Through this Jordan decomposition, the degrees are multiplied by
$|\GG^*{}^{F_0}|_{2'}|\Ce G(s)|^{-1}_{2'}$, so
$|\Irr_{2'}(G)|=\sum_s |\ser{\Ce G(s)}1_{2'}|$, a sum over the semi-simple
classes of $G$.

Characteristic polynomials provide a bijection between the classes of semi-simple elements of $\Sp_{2n}(\F_2)$ and the set of self dual polynomials $f\in\F_2[X]$ of degree $2n$. If
$s$ corresponds with $f$, then
$\Ce G (s)\cong \Sp_{2m}(\F_2)\times C_s$ where $C_s$ is a product
of finite linear groups and $2m$ is the multiplicity of $(X-1)$ in $f$.
For a given $m<n$, the number of such classes is $2^{n-m-1}$. This
is because one has to count the polynomials
$f=(X-1)^{2m}g$ with a self dual $g(X)=1+a_1X+\dots +
a_{n-m-1}X^{n-m-1}+a_{n-m}X^{n-m}+a_{n-m-1}X^{n-m+1} +\dots
+a_{1}X^{2n-2m-1}+X^{2n-2m}$ such that $g(1)\not=0$. Such $g$'s
are $2^{n-m-1}$, corresponding to the choice of coefficients at
degrees $1,2,\dots , n-m-1$ since
$g(1)=a_{n-m}$ has to be $=1$. For $m=n$ (central element) there
is 1 conjugacy class ($s=1$).

The unipotent characters of finite reductive groups of type $A$ are
of even degrees except the trivial character (see for
instance [H] or [M1]~6.8). Then Proposition~1 implies that each semi-simple
class
$s$ corresponding with $m$ as above satisfies
$|\ser{\Ce G(s)}1_{2'}|=5$ for $m\geq 2$, $|\ser{\Ce G(s)}1_{2'}|=1$
otherwise. So the above implies that
$$|\Irr_{2'}(G)|=5.\sum_{m=2}^{n-1}2^{n-m-1}+5+2^{n-2}+2^{n-1}=5.2^{n-2}+3.2^{n-2}=2^{n+1}.$$
This is our claim.
\qed

\medskip

\subsection{The local case.}

We use the
description of Sp$_{2n}(\F_2 )\subset {\rm GL}_{2n}(\F_2)$ as the
subgroup of matrices
$u$ such that
$^tu\left(\begin{matrix} 0&J\\  J&0 \end{matrix} \right)u=\left(\begin{matrix} 0&J\\
J&0 \end{matrix} \right)$ where $J$ denotes the matrix with coefficients $(\dd_{i, n+1-j})_{1\leq i,j\leq n}$ and $u\mapsto {}^tu$ denotes transposition (see [DM] 15.2).

Let $U:=\Big\{
\left(\begin{matrix}  x&xsJ\\  0&J\bar x J\end{matrix} \right)\mid x\in V\ ,\  s\in {\rm
Sym}_n\Big\}$ where ${\rm Sym}_n$ (resp. $V$) is the set of
symmetric (resp. upper triangular unipotent) matrices of order $n$
with coefficients  in
$\F_2$, and one denotes $\bar x = {}^tx^{-1}$. We have

\medskip\noindent {\bf Proposition~3. } {\sl $U$ is a Sylow
$2$-subgroup of
$G=\Sp_{2n}(\F_2)$ for $n\geq 2$. Moreover ${\rm N}_G(U)=U$ and $U/[U,U]$ is of
order
$2^{n+1}$.}
\medskip

\medskip\noindent {\bf Corollary~4. } {\sl Mac Kay conjecture (on
character degrees) is satisfied in $G=\Sp_{2n}(\F_2)$ for the prime 2 ($n\geq 2$). That
is, the
normalizer of any Sylow $2$-subgroup of $G$ has the same number
of characters of odd degrees as $G$ itself.}

\medskip\noindent {\it Proof. } By Proposition~3, the irreducible
characters of
N$_G(U)=U$ of odd degrees are exactly the linear characters of $U$. So their number
is the cardinality of $U/[U,U]$, that is $2^{n+1}$ thanks to
Proposition~3 again. Combining with Proposition~2 gives our claim.
\qed

\medskip\noindent {\it Proof of Proposition~3. }  Note that $U$ equals the group of
elements over $\F_2$
of a rational Borel subgroup (see [DM] 15.2), so it equals its normalizer by the axioms of finite
BN-pairs which are satisfied by this group. Thus our first claim.

Note also the semi-direct decomposition $U\cong {\rm Sym}_n\rtimes
V$ for the action of  $ V$ on ${\rm Sym}_n$ given by
$x.s = xs{}^t\! x$ for $x\in V$, $s\in {\rm Sym}_n$.  Since ${\rm Sym}_n$ is abelian
and since the Sylow
$2$-subgroup
$V$ of GL$_n(\F_2)$ is known to satisfy $|V/[V,V]|=2^{n-1}$ (see for
instance [DM] p.129 and [H]), our claim about $U/[U,U]$ reduces to show that
${\rm Sym}_n/[{\rm Sym}_n,V]$ if of order 4.
So we have to prove that the sum
$S'=\sum_{x\in V}\th_x({\rm Sym}_n)$ of images of endomorphisms
$\th_x\colon s\mapsto xs{}^t\! x-s$ of ${\rm Sym}_n$ has
codimension 2.

For $1\leq i,j\leq n$, let us denote by $E_{ij}$ the usual elementary
matrix of order $n$. We have $E_{ij}+E_{ji}+E_{ii}\in S'$ for any $1\leq
i<j\leq n$, by computing $\th_x(s)$ for $s=E_{jj}$,
$x=I_n+E_{ij}$. We also have $E_{ij}+E_{ji}\in S'$ for any $1\leq
i<j\leq n$ with $(i,j)\not= (n-1,n)$ (taking
$s=E_{jk}+E_{kj}$ and $x=I_n+E_{ik}$ for some $k>i$, $k\not= j$).
This shows that $S'$ contains the $E_{ij}+E_{ji}$'s for $1\leq i<j\leq n$
with $(i,j)\not= (n-1,n)$, along with $E_{11},E_{22},\dots
E_{n-2,n-2}$ and $E_{n-1,n}+E_{n,n-1}+E_{n-1,n-1}$. This makes a
subspace of codimension 2 in ${\rm Sym}_n$, a supplement
subspace being generated by $E_{n-1,n-1}$ and $E_{n,n}$. The action
of
$V$ on the quotient is easily checked to be trivial (one just has to
check the images of
$E_{n-1,n-1}$ and $E_{n,n}$ by $\th_x$ for $x=I_n+E_{ij}$ -- which we
just did above -- since the latter generate
$V$ as a group, using again the fact that the field has two
elements). So this subspace is indeed the sum of the images of all
the $\th_x$'s for $x\in V$.
\qed

\medskip

\subsection { Inductive McKay condition for
{\rm Sp}$_{2n}(2)$}

In [IMN], the authors show that a finite group satisifies the McKay conjecture as soon as
all its simple non-abelian subquotients satisfy a series of conditions concerning the
automorphism groups of the perfect central extensions of those simple subquotients. Note that $\Sp_{4}(\F_2)$ is isomorphic with the symmetric group on 6 letters.

\medskip\noindent {\bf Theorem~5. } {\sl Let
$n\geq 3$ be an integer. Then $\Sp_{2n}(\F_2)$ is a simple group that satisfies the
conditions of [IMN]~\S~10 for all prime numbers.}

\medskip

\medskip\noindent {\it Proof of the Theorem. } When $n\geq 3$, $\Sp_{2n}(\F_2)$ is a
simple group. When $n=3$, it satisfies the theorem by [M2]~4.1. When $n>3$, 
$\Sp_{2n}(\F_2)$ has trivial Schur multiplier and trivial outer automorphism group (see 
[GLS]), so that the checking required by [IMN] just amounts to the McKay conjecture
itself (see [IMN]~10.3). For $\ell =2$, it is Corollary~4. In the case of other
primes, this is a consequence of Malle's parametrization [M1]~7.8 along with Sp\"ath's extensibility
results (see [S1]~1.2, [S2]~1.2, 8.4).
\qed

\bigskip

\section{$\Sp_4(\F_{2^m})$.}

We prove here the following

\medskip\noindent {\bf Theorem~6. } {\sl Let
$m\geq 2$ be an integer. Then $\Sp_{4}(\F_{2^m})$ is a simple group that satisfies the
conditions of [IMN]~\S~10 for all prime numbers.}

\medskip

We keep $\F$ the algebraic closure of the field $\F_2$ with 2 elements and $\GG =\Sp_4(\F )$. We denote by $\TT_0$ its diagonal torus and $(\F ,+)\to\GG$, $t\mapsto x_\aa (t)$ its minimal unipotent $\TT_0$-stable subgroups indexed by the $\TT_0$-roots $\aa$. The Weyl group $W:=\No\GG (\TT_0)/\TT_0$ is generated by the classes $s_1$, $s_2$ of the permutation matrices in $\GL_4(\F )$ associated with the permutations $(1,2)(3,4)$ and $(2,3)$, respectively.

Denote by $F_0'$ the automorphism of $\GG =\Sp_4(\F )$ which sends $x_\aa (t)$ to $x_{\aa '}(t^2)$ if $\aa $ is short (i.e. its associated reflection is conjugated with $s_1$), to $x_{\aa '}(t)$ otherwise, and where $\aa\mapsto \aa '$ is the permutation of roots corresponding to the swap of $s_1$ and $s_2$, see [C1]~12.3.3. Note that $F_0=(F'_0)^2$ (notation of \S 1). Denote $F=F_0^m$, so that $\Sp_{4}(\F_{2^m})=\GG^F$.

\medskip

\medskip\noindent\it Proof of Theorem~6. \rm The group $G=\Sp_4(\F_{2^m})$ ($m\geq 2$) is simple with trivial Schur multiplier and cyclic outer automorphism group generated by $F'_0$ (see [GLS]). Then the conditions of [IMN]~\S~10 amount to find for each prime $\ell$ dividing $|G|$ a proper subgroup $N<G$ containing $\No G(P)$ for $P$ a Sylow $\ell$-subgroup of $G$ and such that $\si (N)=N$ and $|\Irr_\lp (G)^\si |=|\Irr_\lp (N)^\si|$ for any $\si\in\No{{\Aut(G)}}(P)$ (see [Br]~\S~3). The case of $\ell =2$ is also done in [Br], so we assume that $\ell$ is odd dividing $(2^{4m}-1)(2^{2m}-1)=|\Sp_4(\F_{2^m})|_{2'}$. The order of $2^m$ mod. $\ell$, is $e\in\{ 1,2,4\}$. Let $\SS_e$ be a Sylow $\phi_e$-torus of $\GG$. We have that $\TT_e:=\Ce\GG (\SS_e )$ is a maximal torus of $\TT_0$-type $w_e=$ 1, $s_1s_2s_1s_2$, or $s_1s_2$ according to $e$ being 1, 2 or 4 (for types of maximal $F$-stable tori, and latter Levi subgroups, we refer to [DM] p.113).

Arguing as in the proof of [M1]~5.14, any Sylow $\ell$-subgroup $P$ has a unique maximal toral elementary abelian subgroup whose normalizer $N$ in $G$ is then also $N:=\No G(\SS_e )=\No G(\TT_e )$. It is stable by any automorphism $\si$ such that $\si (P)=P$. From what has been said about possible $\si$'s, and noting that $N$ has an abelian normal subgroup $\TT_e^F$ with $\lp$ index, we see that we must just prove that $$|\Irr_\lp (G)^{F'}|=|\Irr (N)^{xF'}|\leqno({\rm E})$$ for any $F'$ a power of $F'_0$ and some $x\in G$ is such that $F'(\SS_e )=\SS_e^x$.

Bringing $(\TT_e, F)$ to $(\TT_0 ,w_eF)$ by conjugacy with some $g\in\GG$ such that $g^{-1}F(g)\in w_e\TT_0$, we may rewrite the above as $$|\Irr_\lp (G)^{F''}|=|\Irr (\No\GG(\TT_0)^{w_eF})^{F''}|\leqno({\rm E'})$$ if $F''$ commutes with $w_eF$ and is in the same class as $F'$ mod inner automorphisms of $G$.

Recall Malle's bijection $\Irr_\lp (G)\mr{\sim}\Irr_\lp (N)$ which, among other properties, sends components of $\Lu{\TT_e}\GG\th$ to components of $\Ind_{\TT_e}^{N}\th$ for relevant $\th\in\Irr (\TT_e^F)$ (see [M1]~\S 7.1).

Let us first look at regular characters $\pm\Lu\TT\GG (\th )$. They are of degree $\lp$ if and only if $\TT$ can be taken as $\TT_e=\Ce\GG (\SS_e )$ (see [M1]~6.6). Such a character is fixed by $F'$ if and only if $F'(\TT_e ,\th )$ and $(\TT_e ,\th )$ are $\GF$-conjugate (see [Br]\S~2.1.2). This is equivalent to $xF'(\th )$ being $\No G(\SS_e)$-conjugate to $\th$ ([M1]~5.11). This is also the criterion for $\Ind_\TF^N (\th )$ being $xF'$-fixed as can be seen easily from the definition of induced characters. Thus our claim in the form of (E) above.

Let us now turn to unipotent characters. From [M1] 6.5, we know that they have to be in $\ser\GF {\TT_e ,1}$, the set of irreducible characters occuring in the generalized character $\Lu{\TT_e}\GG 1$. So we have to check that $\ser\GF {\TT_e ,1}^{F'}_\lp$ and $\Irr(N/\TT_e^F )^{F'}$ have same cardinality.

As for the first set, one knows that among the six unipotent characters of $\Sp_4(\F_{2^m})$, only the two that are of generic degree ${1\over 2}q(q^2+1)$ are not fixed by $F'_0$ (see [M1]~3.9.a). Those are among unipotent characters of degree prime to $\ell$ only when $e=1$ or $2$. So it suffices to check that all characters of $N/\TT_e^F $ but 2 are fixed by $xF'$ in case $e=1$ or $2$ and $F'$ is an odd power of $F'_0$, and that all are fixed otherwise.

In cases $e=1$ or $2$,  $w_1=1$, $w_2=s_1s_2s_1s_2$ both are fixed by $F'_0$, so one may take $F''=F'$ in (E') above. Recall that $F_0'$ acts on $W$ by permuting $s_1$ and $s_2$. The group $W$ is dihedral of order 8, so $F'_0$ induces an automorphism of order two of $W^{\rm ab}$, so two linear characters out of four are $F_0'$-fixed, while the character of degree two is fixed. Hence our claim for any odd power of $F'_0$. In the case of an even power, the action is trivial, as expected.

In the case $e=4$, one may take $w_4=s_1s_2$ and $F''=(s_1F'_0)^a$ when $F' =(F'_0)^a$. Then the action of $F''$ on $(\No\GG (\TT_0 )/\TT_0)^{w_4F} =\Ce W(w_4)$ is trivial.

We now assume $\ser Gs^{F'}_\lp\not=\vid$ for an $s$ that is neither central nor regular. 
The group $\Ce\GG (s)$ is always a Levi subgroup of $\GG$ (see proof of Proposition~2 above) and by [M1] 6.5 it must contain a Sylow $\phi_e$-torus. A proper $F$-stable Levi subgroup of $\GG$ can contain a $\phi_1$-Sylow for types $(\LL_{\{s_1\}}, F)$ and $(\LL_{\{s_2\}},  F)$ and a $\phi_2$-Sylow for types $(\LL_{\{s_1\}}, s_2s_1s_2F)$ and $(\LL_{\{s_2\}}, s_1s_2s_1F)$. In each case the corresponding finite group has two unipotent characters, the trivial and the Steinberg characters, of distinct degrees, so that for an $s$ whose class is $F'$-stable with such a centralizer in the dual, $\ser Gs$ has two elements with distinct degrees, so $F'$ acts trivially on $\ser Gs$. 

The corresponding statement on the local side is as follows : if $\th$ is a non regular non central linear character of $\TT_0^{w_eF}$, then $\Ind^{\No\GG(\TT_0)^{w_eF}}_{\TT_0^{w_eF}} \th$ has two elements both $F''$-fixed if $F''(\th) \in \No\GG(\TT_0)^{w_eF}.\th$. This holds because non-regularity implies 
$(\No\GG(\TT_0)^{w_eF})_\th/\TT_0^{w_eF}$ is of order 2, but then $F''$ can act only trivially on it.
\qed

\medskip

\section{References.}
\begin{enumerate}
%\item[{[Be]}] D. Benson, \sl Representations and Cohomology I : Basic representation Theory of Finite Groups and Associative Algebras, \rm Cambridge, 1997.

\item[{[Br]}] O. Brunat, On the inductive McKay condition in the defining
characteristic, {\sl preprint}.

\item[{[CS]}] M. Cabanes and B. Sp\"ath, Equivariance and extendibility in finite reductive groups with connected center, \it in preparation, \rm 2011.

\item[{[C1]}] R. Carter, {\it Simple groups of Lie type}, Wiley, New York 1972.

\item[{[C2]}] R. Carter, {\it Finite groups of Lie type : conjugacy classes and
complex characters}, Wiley, New York 1985.

\item[{[DM]}] F. Digne and J. Michel, {\it Representations of finite groups of
Lie type}, Cambridge, 1991.

\item[{[GLS]}]  D. Gorenstein, R. Lyons and R. Solomon, {\it The classification of the finite
simple groups,
Number 3.} {\sl Mathematical Surveys and Monographs, Amer. Math. Soc.}, Providence, 1998.

\item[{[H]}] R. Howlett, On the degrees of Steinberg characters of Chevalley
groups, {\sl Math. Z.}, {\bf 135} (1974), 125-135.

\item[{[IMN]}] M. Isaacs, G. Malle, and G. Navarro, A reduction theorem for
McKay conjecture, \sl Inventiones Math., \bf 170\rm (2007),
33--101.

\item[{[L]}] G. Lusztig, Irreducible representations of finite classical
groups, \sl Inventiones Math.
\bf 43 \rm (1977), 125--175.

\item[{[M1]}]  G. Malle, Height 0 characters of finite groups of Lie type, \sl
Representation Theory,
\bf 11\rm (2007), 192--220.

\item[{[M2]}]  G. Malle, The inductive
McKay condition for simple groups not of Lie type, \sl Comm.
Algebra.,
\bf 36-2\rm (2008), 455--463.

%\item[{[M3]}] ??  G. Malle, Extensions of unipotent characters and the inductive McKay condition, \sl J. Algebra, \bf 320\rm (2008), 2963-2980.

\item[{[S1]}] B. Sp\"ath, Sylow $d$-tori of classical groups and the McKay conjecture I, \sl
J. Algebra, \bf 323\rm (2010), 2469-2493.

\item[{[S2]}] B. Sp\"ath, Sylow $d$-tori of classical groups and the McKay conjecture II,\sl
J. Algebra, \bf 323\rm (2010), 2494-2509.

\item[{[S3]}] B. Sp\"ath, Inductive McKay condition in defining characteristic, \it preprint \rm (2010) arXiv:1009.0463 .

\end{enumerate}

\end{document}